\documentclass[11pt,reqno]{amsart}

\usepackage{amsmath,amsthm}
\usepackage{amssymb}
\usepackage{euscript}
\usepackage{url}

\numberwithin{equation}{section}

\newtheorem{Thm}[subsection]{Theorem}
\newtheorem{Lem}[subsection]{Lemma}

\newtheorem{Cor}[subsection]{Corollary}

\theoremstyle{definition}

\begin{document}
\title[]{Characteristic rank of vector bundles over Stiefel manifolds}
\author{J\'ulius Korba\v{s}, Aniruddha C. Naolekar, and Ajay Singh Thakur}
\thanks{Part of this research was carried out while J. Korba\v s was a
member of two research teams supported in part by the grant agency VEGA (Slovakia).}


\begin{abstract}
The characteristic rank of a vector bundle $\xi$ over a finite connected $CW$-complex $X$ is
by definition the largest integer $k$, $0\leq k\leq \mathrm{dim}(X)$, such that every
cohomology class $x\in H^j(X;\mathbb Z_2)$, $0\leq j\leq k$, is a polynomial in the
Stiefel-Whitney classes $w_i(\xi)$. In  this note we compute the characteristic rank of
vector bundles over the Stiefel manifold $V_k(\mathbb F^n)$, $\mathbb F=\mathbb R,\mathbb
C,\mathbb H$.

\end{abstract}


\email{}
\address{J. K., Department of Algebra, Geometry, and Mathematical Education,
Faculty of Mathematics, Physics, and Informatics, Comenius University, Mlynsk\'a dolina,
SK-842 48 Bratislava 4, SLOVAKIA or Mathematical Institute, Slovak Academy of Sciences,
\v{S}tef\'anikova 49, SK-814 73 Bratislava 1, SLOVAKIA} \email{ korbas@fmph.uniba.sk}

\address{A. C. N. and A. S. T., Indian Statistical Institute,
8th Mile, Mysore Road, RVCE Post, Bangalore 560059, INDIA.}
\email{ani@isibang.ac.in, thakur@isibang.ac.in}

\date{}
\maketitle

\section{Introduction}

Let $X$ be a connected finite $CW$-complex and $\xi$ a real vector bundle over $X$. Recall
\cite{nt} that the {\em characteristic rank of} $\xi$ over $X$, denoted by
$\mathrm{charrank}_{X}(\xi)$, is by definition the largest integer $k$, $0\leq k\leq
\mathrm{dim}(X)$, such that every cohomology class $x\in H^j(X;\mathbb Z_2)$, $0\leq j\leq
k$, is a polynomial in the Stiefel-Whitney classes $w_i(\xi)$. The {\em upper characteristic
rank of} $X$, denoted by $\mathrm{ucharrank}(X)$, is the maximum of
$\mathrm{charrank}_{X}(\xi)$ as $\xi$ varies over all vector bundles over $X$.

Note that if $X$ and $Y$ are homotopically equivalent connected closed manifolds, then
$\mathrm{ucharrank}(X)=\mathrm{ucharrank}(Y)$. When $X$ is a connected closed smooth manifold
and $TX$ the tangent bundle of $X$, then $\mathrm{charrank}_{X}(TX)$, denoted by
$\mathrm{charrank}(X),$ is called the {\em characteristic rank of the manifold} $X$ (see
\cite{korbas}).

The characteristic rank of vector bundles can be used to obtain bounds for the $\mathbb
Z_2$-cup-length of manifolds (see  \cite{balkokorbas},\cite{korbas} and \cite{nt}). In some
situations, the value of the upper characteristic rank can be used to show the vanishing of
the Stiefel-Whitney class of a certain degree for all vector bundles. An important task is
therefore to understand the characteristic rank of vector bundles.

In \cite{nt}, the second and third named authors have computed the characteristic rank of
vector bundles over: a product of spheres, the real  and complex projective spaces, the Dold
manifold $P(m,n)$, the Moore space $M(\mathbb Z_2,n)$, and the stunted projective space
$\mathbb R\mathbb P^n/\mathbb R\mathbb P^m$.

Let $\mathbb F$ denote either the field $\mathbb R$ of reals, the field $\mathbb C$ of
complex numbers or the skew-field $\mathbb H$ of quaternions. Let $V_k(\mathbb F^n)$ denote
the Stiefel manifold of orthonormal $k$-frames in $\mathbb F^n$. In this note we compute the
characteristic rank of vector bundles over $V_k(\mathbb F^n)$. Our methods are elementary and
make use of some well-known facts about the Stiefel manifolds. We prove the following.

\begin{Thm}\label{maintheorem}
Let $X=V_k(\mathbb F^n)$ with $1<k<n$ when $\mathbb F=\mathbb R$ and
$1<k\leq n $ when $\mathbb F=\mathbb C, \mathbb H$.
\begin{enumerate}
\item If $\mathbb F=\mathbb R$, then
$$\mathrm{ucharrank}(X)=\left\{\begin{array}{cl}
n-k-1 & \mbox{if $n-k\neq 1,2,4,8$,}\\
2 & \mbox{if $n-k=1$ and $n\geq 4$,}\\
2 & \mbox{if $n-k=2$,}\\
4 & \mbox{if $n-k=4$ and $k=2$.}\end{array}\right.$$
\item If $\mathbb F=\mathbb R$, $k>2$ and $n-k=4$, then $\mathrm{ucharrank}(X)\leq 4$.
\item If $\mathbb F=\mathbb R$ and $n-k=8$, then $\mathrm{ucharrank}(X)\leq 8$.
\item If $\mathbb F=\mathbb C$, then
$$\mathrm{ucharrank}(X)=\left\{\begin{array}{cl}
2 & \mbox{if $k=n$,}\\
2(n-k) &\mbox{if $k<n$.}\end{array}\right.$$
\item If $\mathbb F=\mathbb H$, then $\mathrm{ucharrank}(X)=4(n-k)+2$.
\end{enumerate}
\end{Thm}

When $\mathbb F=\mathbb R$ and $n-k=4,8$, we only give a bound.

The characteristic rank of vector bundles over $V_1(\mathbb F^n)$, which
is a sphere, and $SO(3)=\mathbb R\mathbb P^3$ has been described in \cite{nt}.
Note that $V_n(\mathbb R^n)=O(n)$ is not connected. This explains the restrictions
on the range of $k$ and the condition $n\geq 4$ in the second
equality of (1) in the above theorem.

Notations. The characteristic rank of a vector bundle $\xi$ over $X$ will simply be denoted
by $\mathrm{charrank}(\xi)$; the space $X$ will usually be clear from the context. For a
space $X$, $H^*(X)$ will denote cohomology with $\mathbb Z_2$-coefficients.

\section{Proof of Theorem\,\ref{maintheorem}}

We begin by recalling certain standard facts about the Stiefel manifolds that are needed to
prove the main theorem. One fact about Stiefel manifolds that we shall need is a description
of the $\mathbb Z_2$-cohomology ring of $V_k(\mathbb F^n)$. We note this below.

\begin{Thm}\label{borel}{\rm (\cite{borel}, Propositions 9.1 and 10.3)} We have
$H^i(V_k(\mathbb F^n))=0$ for $i=1,2,\ldots , c(n-k+1)-2$ and $H^{c(n-k+1)-1}(V_k(\mathbb
F^n))\cong \mathbb Z_2$, where $c=\dim _{\mathbb R}\mathbb F$. Further, when $\mathbb
F=\mathbb R$, the cohomology ring $H^*(V_k(\mathbb R^n))$ has a simple system of generators
$a_{n-k}, a_{n-k+1},\ldots , a_{n-1}$ $(a_i\in H^i(V_k(\mathbb R^n)))$ such that
$a_i^2=a_{2i}$ if $2i\leq n-1$ and $a_i^2=0$ otherwise. \qed
\end{Thm}

The action of the Steenrod squares on $H^*(V_k(\mathbb R^n)$ is given
by (see \cite{borel}, Remarque 2 in \S 10)
$$Sq^i(a_j)=\left\{\begin{array}{cl}
\displaystyle{j \choose i}a_{j+i} & \mbox{if $ j+i\leq n-1$,}\\
& \\
0 & \mbox{otherwise.}
\end{array}\right.$$

For $k\geq 2$ consider the sphere bundle $S^{n-k}\stackrel{i}\hookrightarrow V_k(\mathbb
R^n)\stackrel{p}\longrightarrow V_{k-1}(\mathbb R^n)$, where $p$ maps a $k$-frame to the
$(k-1)$-frame determined by ignoring the last vector. It is clear, by the above theorem, that
the Serre spectral sequence of this sphere bundle is trivial and hence, the homomorphism
$i^*:H^{n-k}(V_k(\mathbb R^n))\longrightarrow H^{n-k}(S^{n-k})$ is an isomorphism.

\begin{Lem}\label{lemma} Provided that $k\neq n$ for $\mathbb F = \mathbb C$,  let
$n$ and $k$ be as in Theorem\,\ref{maintheorem}. Then, for any vector bundle $\xi$ over
$V_k(\mathbb F^n)$, we have
\begin{enumerate}
\item $\mathrm{charrank}(\xi)\geq c(n-k+1)-2$,
\item $\mathrm{charrank}(\xi)=c(n-k+1)-2$, if $\mathbb F$ is $\mathbb C$ or $\mathbb H$,
\item $\mathrm{charrank}(\xi)\leq n-k$ if $\mathbb F$ is $\mathbb R$ and $n-k$ is even.

\end{enumerate}
\end{Lem}
{\bf Proof.} (1) follows from the cohomology structure of $V_k(\mathbb F^n)$. Next observe that by Wu's formula
$$w_{c(n-k+1)-1}(\xi)=w_1(\xi)w_{c(n-k+1)-2}(\xi) + Sq^1(w_{c(n-k+1)-2}(\xi))=0.$$
This proves (2). We now come to the proof of (3). Note that, by
Theorem\,\ref{borel}, $H^{n-k}(V_k(\mathbb R^n)))\cong \mathbb Z_2$ is generated by $a_{n-k}$ and
$H^{n-k+1}(V_k(\mathbb R^n))\cong \mathbb Z_2$ is generated by $a_{n-k+1}$. Suppose that $\mathrm{charrank}(\xi)\geq n-k+1$. Then
$w_{n-k}(\xi)=a_{n-k}$ and $w_{n-k+1}(\xi)=a_{n-k+1}$. Now by Wu's formula we have
$$w_{n-k+1}(\xi)=w_1(\xi)w_{n-k}(\xi)+Sq^1(w_{n-k}(\xi))=(n-k)a_{n-k+1}=0.$$
This contradiction proves (3).\qed

We are now in a position to prove our main theorem.

{\em Proof of Theorem\,\ref{maintheorem}.} If $\xi$ is a vector bundle over $V_k(\mathbb
R^n)$ with $w_{n-k}(\xi)\neq 0$, then $i^*\xi$ is a vector bundle over $S^{n-k}$ with
$w_{n-k}(i^*\xi)\neq 0$. By Theorem 1 in \cite{milnor}, this is possible only if $n-k=
1,2,4,8$. This and Lemma \ref{lemma}(1) prove the first equality in Theorem
\ref{maintheorem}(1).

To prove the second equality in Theorem\,\ref{maintheorem}(1), we note that
$H^1(V_{n-1}(\mathbb R^n))=H^1(SO(n))$ is generated by $a_1$, $H^2(SO(n))$ is generated by
$a_2=a_1^2$ and $H^3(SO(n))$ is generated by $a_1^3=a_1a_2$ and $a_3$. Now if $\xi$ is a
non-orientable line bundle over $SO(n)$, then clearly $\mathrm{charrank}(\xi)\geq 2$. Now
assume that $\xi$ is a vector bundle over $SO(n)$ with $\mathrm{charrrank}(\xi)\geq 3$. Then
$w_1(\xi)=a_1$, $w_2(\xi)=ka_2=ka_1^2$ with $k\in\{0,1\}$ and $w_3(\xi)=a_3$ is not a
multiple of $a_1^3$. But by Wu's formula
$$w_3(\xi)=w_1(\xi)w_2(\xi)+Sq^1(w_2(\xi))=a_1^3.$$
This contradiction completes the proof.

To prove the third equality in Theorem\,\ref{maintheorem}(1), we note that by
Lemma\,\ref{lemma}(3), we have $\mathrm{charrank}(\xi)\leq 2$ for any vector bundle $\xi$
over $V_k(\mathbb R^n)$ when $n-k=2$. It is well known that $H^2(V_{n-2}(\mathbb R^n);\mathbb
Z)\cong \mathbb Z$. Then clearly there is a $2$-plane bundle $\xi$ over $V_{n-2}(\mathbb
R^n)$ with Euler class $e(\xi)$ a generator and hence $w_2(\xi)\neq 0$. This implies that
$\mathrm{charrank}(\xi)=2$. This completes the proof.

In view of Lemma\,\ref{lemma}(3), the proof of the fourth equality in
Theorem\,\ref{maintheorem}(1) will be complete if we exhibit a vector bundle $\xi$ over
$V_2(\mathbb R^6)$ with $w_4(\xi)\neq 0$. To construct such a bundle, we start with the
well-known circle bundle $p:V_2(\mathbb R^6)\longrightarrow G_2(\mathbb C^4)$. Here
$G_2(\mathbb C^4)$ denotes the complex Grassmann manifold of complex $2$-planes in $\mathbb
C^4$. Since $H^3(G_2(\mathbb C^4))=0$, the Gysin sequence
$$\cdots\rightarrow H^2(G_2(\mathbb C^4))\stackrel{\psi}\longrightarrow H^4(G_2(\mathbb C^4))
\stackrel{p^*}\longrightarrow H^4(V_2(\mathbb R^6))\longrightarrow H^3(G_2(\mathbb
C^4))\rightarrow\cdots$$ of the circle bundle $p$ shows that the homomorphism
$p^*:H^4(G_2(\mathbb C^4))\longrightarrow H^4(V_2(\mathbb R^6))$ is onto. Let $\gamma$ be the
canonical complex $2$-plane bundle over $G_2(\mathbb C^4)$ and $\gamma_{\mathbb R}$ its
underlying real bundle. It is known that $H^2(G_2(\mathbb C^4))\cong \mathbb Z_2$ is
generated by $w_2(\gamma_{\mathbb R})$ and $H^4(G_2(\mathbb C^4))\cong \mathbb Z_2\oplus
\mathbb Z_2$ is generated by $w_2^2(\gamma_{\mathbb R})$ and $w_4(\gamma_{\mathbb R})$. Since
$\psi(x)=x\smile (w_2(\gamma_{\mathbb R}))$ and $w_2^2(\gamma_{\mathbb R})\neq
w_4(\gamma_{\mathbb R})$, it is clear that $w_4(\gamma_{\mathbb R})\notin
\mathrm{Im}(\psi)=\mathrm{Ker}(p^*)$, thus $p^*(w_4(\gamma_{\mathbb
R}))=w_4(p^*\gamma_{\mathbb R})\neq 0$. Thus $\xi=p^*\gamma_{\mathbb R}$ is the required
vector bundle.

The assertions in (2) and (3) follow from Lemma\,\ref{lemma}(3).

To prove the assertion (4) in Theorem\,\ref{maintheorem}, first assume that $k=n$ with $n\geq
2$. Then, by Theorem\,\ref{borel}, $H^1(U(n))\cong \mathbb Z_2$. Thus there exists a
non-trivial line bundle $\xi$ such that $w_1(\xi)$ generates $H^1(U(n))$. Since
$H^2(U(n))=0$, it follows that $\mathrm{ucharrank}(U(n))\geq 2$. Now if there exists a vector
bundle $\xi$ with $w_3(\xi)\neq 0$, then $w_3(\xi)$ generates $H^3(U(n))\cong \mathbb Z_2$.
But by Wu's formula we have
$$w_3(\xi)=w_1(\xi)w_2(\xi)+Sq^1(w_2(\xi))=0.$$
This is a contradiction. The case $k<n$ follows from Lemma\,\ref{lemma}(2).

The assertion (5) also follows from Lemma\,\ref{lemma}(2). This completes the proof of the
theorem. \qed

We have the following immediate corollary of Theorem \ref{maintheorem}.

\begin{Cor}\label{corollary}
\begin{enumerate}
\item If $n-k\neq 1,2,4,8$ and $1<k<n$, then we have $w_{n-k}(\xi)=0$ for any vector bundle $\xi$ over $V_k(\mathbb R^n)$.
\item For any non-orientable vector bundle $\xi$ over $SO(n)$, $n\geq 4$, we have either $w_3(\xi)=0$ or
$w_3(\xi)=a_1^3$, where $a_1$ is the $($unique$)$ non-zero element in the first cohomology.
\item Let $\xi$ be a vector bundle over $V_k(\mathbb R^{n})$, where $n-k$ is even and $1 < k < n$. Then $w_{n-k+1}(\xi)=0$.
\item For any vector bundle $\xi$ over $U(n)$, we have $w_3(\xi)=0$.
\end{enumerate} \qed
\end{Cor}

\end{document}